\documentclass{amsart}  
\usepackage{amssymb,amscd}
\begin{document}
 
\title[Integration of Lie algebroids]{Groupoids and the integration of
Lie algebroids}
%
%
\newcommand{\mar}[1]{{\marginpar{\textsf{#1}}}}
\newcommand\datver[1]{\def\datverp%
{\par\boxed{\boxed{\text{Version: #1; Run: \today}}}}}
\datver{3.0A; Revised: 4/29/99}

%
%
\newcommand\Exp{\operatorname{Exp}}
\newcommand\Ran{\operatorname{Im}}
\newcommand\abs{\operatorname{abs}}
\newcommand\rel{\operatorname{rel}}
\newcommand\inv{\operatorname{inv}}
\newcommand\topo{\operatorname{top}}
\newcommand\opp{\operatorname{op}}
\newcommand\mfk{\mathfrak}
\newcommand\coker{\operatorname{coker}}
\newcommand\hotimes{\hat \otimes}
\newcommand\ind{\operatorname{ind}}
\newcommand\End{\operatorname{End}}
\newcommand\per{\operatorname{per}}
\newcommand\pa{\partial}
\newcommand\sign{\operatorname{sign}}
\newcommand\supp{\operatorname{supp}}
\newcommand\cy{\mathcal{C}^\infty}
\newcommand\Cy{\mathcal{C}^\infty}
\newcommand\lra{\longrightarrow}
\newcommand\vlra{-\!\!\!-\!\!\!-\!\!\!\!\longrightarrow}
\newcommand\bS{{}^b\kern-1pt S}
\newcommand\bT{{}^b\kern-1pt T}
\newcommand\Hom{\operatorname{Hom}}

\newcommand\alg[1]{\mathfrak{A}(#1)} 
\newcommand\qalg[2]{\mathfrak{B}_{#1}(#2)} 
\newcommand\ideal[1]{\mathfrak{I}(#1)} 
\newcommand\qideal[2]{\mathfrak{R}_{#1}(#2)} 
\newcommand\In{\operatorname{In}}

\newcommand\TR{\operatorname{T}}
\newcommand\ha{\frac12} 
\newcommand\cal{\mathcal}
\newcommand\END{\operatorname{END}}
\newcommand\ENDG{\END_{\GR}(E)}
\newcommand\CC{\mathbb C}
\newcommand\NN{\mathbb N}
\newcommand\RR{\mathbb R}
\newcommand\ZZ{\mathbb Z}
\newcommand\ci{${\mathcal C}^{\infty}$}
\newcommand\CI{{\mathcal C}^{\infty}}
\newcommand\CIc{{\mathcal C}^{\infty}_{\text{c}}}
\newcommand\hden{{\Omega^{\lambda}_d}}
\newcommand\VD{{\mathcal D}}
\newcommand\mhden{{\Omega^{-1/2}_d}}
\newcommand\ehden{r^*(E)\otimes {\Omega^{\lambda}_d}}


\newcommand{\Cat}{\mathcal C}
\newcommand{\Gr}[1]{{\mathcal G}^{(#1)}}
\newcommand{\GR}{\mathcal G}
\newcommand{\LGR}{\mathcal L}
\newcommand{\GG}{\mathcal G}
\newcommand{\OA}{\mathcal O}
\newcommand{\PS}[1]{\Psi^{#1}(\GR;E)}
\newcommand{\tPS}[1]{\Psi^{#1}(\GR)}
\newcommand{\ttPS}[1]{\Psi_{\loc}^{#1}(\GR;E)}
\newcommand{\tttPS}[1]{\Psi_{\loc}^{#1}(\GR)}
\newcommand{\AL}{{\mathcal A}(\GR)}
\newcommand{\FAM}{P=(P_x,x \in \Gr0)}
\newcommand\symb[2]{{\mathcal S}^{#1}(#2)}
\newcommand{\loc}{\operatorname{loc}}
\newcommand{\cl}{\operatorname{cl}}
\newcommand{\A}{s}
\newcommand{\prop}{\operatorname{prop}}
\newcommand{\comp}{\operatorname{comp}}
\newcommand{\adb}{\operatorname{adb}}
\newcommand{\dist}{\operatorname{dist}}

\newcommand{\alp }{r }
\newcommand{\bet }{d }
\newcommand{\gm }{\Gamma }
\newcommand{\lon }{\longrightarrow }
\newcommand{\be }{\begin{eqnarray*}}
\newcommand{\ee }{\end{eqnarray*}}
\newcommand{\GGR}{{\GR}}
\newcommand{\cald}{{\cal D}}
\newcommand{\calx}{{\cal X}}

\def\nin{\noindent}
\def\eg{e.g.\ }
\def\pt#1#2{{\partial #1\over \partial #2}}

\newcommand{\frakg}{{\mathfrak g}}

\let\Tilde=\widetilde
\let\Bar=\overline
\let\Vec=\overrightarrow
\let\ceV=\overleftarrow
\def\vlra{\hbox{$\,-\!\!\!-\!\!\!-\!\!\!-\!\!\!-\!\!\!
-\!\!\!-\!\!\!-\!\!\!-\!\!\!-\!\!\!\longrightarrow\,$}}

\def\vleq{\hbox{$\,=\!\!\!=\!\!\!=\!\!\!=\!\!\!=\!\!\!
=\!\!\!=\!\!\!=\!\!\!=\!\!\!=\!\!\!=\!\!\!=\!\!\!=\!\!\!=\,$}}

\def\lrah{\hbox{$\,-\!\!\!-\!\!\!
-\!\!\!-\!\!\!-\!\!\!-\!\!\!-\!\!\!\longrightarrow\,$}}

\def\surj{-\!\!\!-\!\!\!-\!\!\!\gg}

\def\inj{>\!\!\!-\!\!\!-\!\!\!-\!\!\!>}


 %
%
\newcommand\Mand{\text{ and }}
\newcommand\Mandset{\text{ and set}}
\newcommand\Mas{\text{ as }}
\newcommand\Mat{\text{ at }}
\newcommand\Mfor{\text{ for }}
\newcommand\Mif{\text{ if }}
\newcommand\Min{\text{ in }}
\newcommand\Mon{\text{ on }}
\newcommand\Motwi{\text{ otherwise }}
\newcommand\Mwith{\text{ with }}
\newcommand\Mwhere{\text{ where }} 
\newcommand\ie{{\em i.e., }} 
%
%
\newtheorem{theorem}{Theorem}
\newtheorem{proposition}{Proposition}
\newtheorem{corollary}{Corollary}
\newtheorem{lemma}{Lemma}
\newtheorem{definition}{Definition}
\theoremstyle{remark}
\newtheorem{remark}{Remark}

%
%

\author[V. Nistor]{Victor Nistor} \address{Pennsylvania State
        University, Math. Dept., University Park, PA 16802}
        \email{nistor@math.psu.edu} \thanks{Partially supported an NSF
        Young Investigator Award DMS-9457859 and a Sloan Research
        Fellowship.
} 
 
\dedicatory\datverp
\begin{abstract} 
We show that a Lie algebroid on a stratified manifold is integrable
if, and only if, its restriction to each strata is integrable. These
results allow us to construct a large class of algebras of
pseudodifferential operators. They are also relevant for the
definition of the graph of certain singular foliations of manifolds
with corners and the construction of natural algebras of
pseudodifferential operators on a given complex algebraic variety.
\end{abstract} 
 
\maketitle 
\tableofcontents

\section*{Introduction\label{Sec.Introduction}}

Differentiable groupoids appear in geometry in various instances, for
example in the theory of connections and parallel transport on fiber
bundles, or in the theory of pseudogroups of transformations.
The concept of groupoid is a generalization of the concept 
of group, the main
difference being that not any two elements of a groupoid are
composable. Intuitively, it is convenient to think of a groupoid as a
set of arrows between various points, called units, two arrows being
composable if, and only if, their ends match. (See Section
\ref{Sec.Basic} for precise definitions.)

One of the main features of differential groupoids is that they are
geometric objects that interpolate between differentiable manifolds
and Lie groups: differentiable manifolds have many units and few
arrows; whereas Lie groups have many arrows and few units--actually only one.
This ``interpolation'' property is valid also at the level of
algebras: to a compact smooth manifold $M$ one associates the
commutative algebra $\CI(M)$ of its differentiable functions; whereas
to a Lie group $G$ one associates the convolution algebra $\CIc(G)$ of
compactly supported smooth functions on the group, which is usually
highly non-commutative. In this way, differentiable groupoids provide a
link between geometry and harmonic analysis. 

The algebras $\CI(M)$ and $\CIc(G)$ are particular cases of the
convolution algebra of a differential groupoid, and this feature makes
groupoids a favorite toy model in non-commutative geometry. From this
point of view, the results of this paper are a first step towards a
generalization of the results of \cite{Brylinski-Nistor1}, from
\'etale groupoids to general differential groupoids.

Recall that a {\em Lie algebroid} is a vector bundle $A \to M$ on a
differentiable manifold $M$ together with a Lie algebra structure on
the space $\Gamma(A)$ of its smooth sections and a Lie algebra
morphism
$$
	\Gamma(A) \rightarrow \Gamma(TM),
$$
defined by a structural vector bundle morphism $q : A \to TM$, called
``the anchor map.''  The Lie algebroid $A$ is called {\em regular}
[resp. {\em transitive}] if, and only if, the anchor map $q:A \to TM$
has locally constant rank [resp. it is onto].  From a classical,
differential geometric point of view, to a differential groupoid there
is associated a ``Lie algebroid,'' which is some sort of
``infinitesimal form'' of the groupoid, generalizing both the Lie
algebra of a group and the tangent space to a differentiable manifold.

In \cite{Douady-Lazard} it was proved that any Lie algebroid such that
$q=0$ is the Lie algebroid of a differentiable groupoid, in other
words, it is {\em integrable}. This result can be regarded as a
generalization of ``Lie's third theorem,'' which states that every
finite dimensional Lie algebra is the Lie algebra of a Lie group.
This is relevant because, unlike (finite dimensional) Lie algebras,
Lie algebroids do not always correspond to differential groupoids,
that is, they are not always integrable. Actually, for a transitive
groupoid, one can define an obstruction to integrability, see
\cite{Mackenzie1}.  Nevertheless, it is still interesting to construct
differentiable groupoids that integrate specific Lie algebroids. Some
examples can be found in \cite{NWX}. In that paper, to a
differentiable groupoid $\GR$ there was associated an algebra of
pseudodifferential operators $\Psi^*(\GR)$, such that, for $\Gr0$
compact, all first order differential operators in $\Psi^*(\GR)$ are
linear combinations of sections of $A(\GR)$ and operators of
multiplication by functions. The integration of Lie algebroids is thus
a first step toward constructing algebras of pseudodifferential
operators, and this explains why we are interested in the problem of
integrating Lie algebroids. See
\cite{connesF,connes,MelroseCongress,Monthubert,NWX,WeinsteinAMS} for more on
the question of constructing pseudodifferential operators.

In this paper, we approach the problem of integrating Lie algebroids
from an abstract point of view, looking for some general methods to
integrate Lie algebroids. Since arbitrary Lie algebroids $\pi: A \to
M$ behave rather wildly, we make two assumptions. First, we assume
that the manifold $M$ has a stratification 
$$
	M =\cup S
$$ 
into disjoint strata, each of which is invariant with respect to the
diffeomorphisms generated by the sections of $A$ and, second, we
assume that the restrictions
$$
	A_S=A\vert_{S},
$$
which are Lie algebroids precisely because the strata are invariant,
are regular. (If $M$ satisfies the first assumption, we say that it
has an ``$A$-invariant regular stratification.'')  

For the class of groupoids with $A$-invariant regular stratifications,
one can approach the ``integration problem'' in two steps. First, a
necessary condition for the integrability of $A$ is the integrability
of each of the restrictions $A_S$, and hence the first step will be to
integrate each of these restrictions (see also below).  Assume then
that we can find a differentiable groupoid $\GR_S$ that integrates
$A_S$, for each $S$. The second step is to ``glue'' the resulting
groupoids $\GR_S$. Surprisingly enough, this na\"{\i}ve approach
actually works in most cases. It works for example if we choose the
integrating groupoids $\GR_S$ to be maximal in a suitable sense
($d$-simply connected), and this is our main general result on the
integration of differential groupoids. As in the paper of Douady and
Lazard \cite{Douady-Lazard}, we obtain in general non-Hausdorff
groupoids.

Since transitive Lie algebroids are a particular case of regular Lie
algebroids, the first part of the problem--that is, integrating
regular Lie algebroids--is similar to the problem of integrating
transitive Lie algebroids, and probably can be handled similarly. In
particular, it is clear that not all regular Lie algebroids are
integrable. It is not our purpose in this paper to study the
integration of general regular algebroids, but we do show how to
integrate particular classes of regular algebroids. For example, we
show that a regular algebroid $A$ is integrable if $\ker(q)$, the
kernel of $q$, consists of semisimple Lie algebras, or if $A$ is a
semi-direct product.

As for the second part of the problem, it turns out that there exists
at most one way to glue the groupoids $\GR_S$ that integrate
$A_S=A\vert_S$, assuming that they exist.  The problem is that the
resulting glued space (a groupoid) is not always a smooth manifold, so
this procedure does not lead directly to a differentiable
groupoid. However, we show that this procedure does lead to a
differentiable groupoid that integrates $A$, provided that all
groupoid $\GR_S$ are $d$-simply connected.  We do not assume that the
strata are regular here. We thus obtain the following result.

\begin{theorem} A Lie algebroid $\pi:A \to M$ on a manifold $M$ with an
$A$-invariant stratification is integrable if, and only if,
it is integrable along each stratum.
\end{theorem}

These results provide us with an explicit way of integrating many Lie
algebroids. As an application of the theorem, we prove the
integrability of certain Lie algebroids on foliated manifolds with
corners. A {\em foliated manifold with corners} is a manifold with
corners, each of whose open faces is a foliated manifold, the
foliations being required to satisfy certain compatibility
relations. This result generalizes a construction due to Winkelnkemper
\cite{Winkelnkemper1}. Previously, Melrose \cite{Melroseb} and
Mazzeo--Melrose \cite{Melrose-Mazzeo} have shown how to integrate
certain particular algebroids. However, their framework was different
from ours.

The results of this paper, together with the results of \cite{NWX},
can be used to construct a natural algebra of pseudodifferential
operators on a complex algebraic variety endowed with a
``$\CI$--resolution of singularities,'' thus making a substantial step
towards a solution of the problem stated in \cite{MelroseCongress}.
Then the methods of \cite{MelroseNistor} can presumably be applied to
study the resulting algebras of operators. Algebras of
pseudodifferential operators on groupoids are also a natural framework
to study adiabatic limits \cite{NWX,Witten1}. The problem of
associating an algebra of pseudodifferential operators to a groupoid
was first formally stated in a paper by Weinstein,
\cite{WeinsteinAMS}. However, before that, in \cite{connesF}, Alain
Connes has constructed algebras of pseudodifferential operators on
foliations, which in our setting corresponds to the case of a regular
groupoid with discrete holonomy (see also \cite{connes}). His methods
have played a role in inspiring the constructions of \cite{NWX}.
 
I would like to thank Alan Weinstein for several useful comments on an
earlier version of this paper. Also, I would like to thank an 
anonymous referee for very carefully reading this paper and for 
several useful suggestions.

\section{Basic concepts\label{Sec.Basic}}

We begin this section by fixing notation and recalling some of the basic
concepts used in this paper.

In the following, we shall use the framework of \cite{NWX}.  In
particular, a {\em manifold} is a smooth manifold, possibly {\em with
corners}. By definition, every point $m$ in a manifold with corners
$M$ has a coordinate neighborhood diffeomorphic to $[0,1)^k \times
\RR^{n-k}$, such that the transition functions are smooth (including
on the boundary). Let $k(m)$ be the least $k$ such that $m$ is in a
set diffeomorphic to $[0,1)^k \times \RR^{n-k}$, and let $\pa_k(M)$ be
the set of points $m$ for which $k(m) = k$. A component of $\pa_k(M)$
is called an open face of codimension $k$. A face of $M$ is the closure of
an open face (of $M$).  A {\em hyperface} of $M$ is a face $H$ of
codimension one. 

It is customary to assume that any hyperface $H$ of a manifold with
corners $M$ is of the form $H=\{x_H = 0 \}$, where $x_H$ is a smooth
positive function on $M$ such that $dx_H \not = 0$ on $H$. If this is
the case, $x_H$ is called a {\em defining function} of $H$. We shall
also assume in this paper that each hyperface of $M$ has a defining
functions, although, most of our results are true even without this
assumption. An {\em interior point} of $M$ is a point of $M$ that
belongs to no hyperface.

A {\em submersion} $f : M \to N$ between two manifolds with corners is
a differentiable map with surjective differential at each point, such
that a non-zero tangent vector to $M$ points inward if, and only if,
its image in $TN$ is non-zero and points inward. By definition, if
$f:M \to N$ is a submersion, then the fibers $f^{-1}(y)$ are smooth
manifolds {\em without} corners. A {\em submanifold with corners} $N
\subset M$ is a closed submanifold $N$ of $M$ such that each face of
$N$ is locally the transverse intersection of $N$ with a face of $M$.

We now define groupoids. Recall that a small category is a category
whose class of morphisms is a set. By definition, a {\em groupoid} is
a small category $\GR$ in which every morphism is invertible.  General
results on groupoids can be found in \cite{Renault1}.

We now fix some notation and make the definition of a groupoid more
explicit. The set of {\em} objects (or {\em units}) of $\GR$ is denoted by
$\GR^{(0)}.$ The set of morphisms (or
{\em arrows}) of $\GR$ is denoted by $\GR^{(1)}=\operatorname{Mor}(\GR).$ We
shall sometimes write $\GR$ instead of $\GR^{(1)}$, by abuse of
notation.  For example, when we consider a space of functions on
$\GR$, we actually mean a space of functions on $\GR^{(1)}$.  We will
denote by $d(g)$ [respectively $r(g)$] the {\em domain} [respectively,
the {\em range}] of the morphism $g:d(g) \to r(g).$ We thus obtain
functions
\begin{equation}
d,r:  \GR^{(1)} \longrightarrow  \GR^{(0)}
\end{equation}
that will play an important role.  The multiplication $\mu: (g, h)
\to\mu(g,h)=g h$ is defined on the set $\GR^{(2)}$ of composable pairs
of arrows:
\begin{equation}\label{Mu.everywhere}
\mu:\GR^{(2)}=\GR^{(1)}\times_M \GR^{(1)}:=\{(g,h):d(g)=r(h)\}
\longrightarrow \Gr1.
\end{equation}
The inversion operation is a bijection $\iota(g)= g^{-1}$ of $\Gr1$.
Denoting by $u(x)$ the identity morphism of the object $x \in \Gr0$,
we obtain an inclusion $u:\Gr0 \to \Gr1$.  We see that a groupoid
$\GR$ is completely determined by the spaces $\Gr0$ and $\Gr1$ and by
the structural morphisms $d,r,\mu,u$, and $\iota$.  We sometimes write
$$
	\GR=(\Gr0,\Gr1,d,r,\mu,u,\iota).
$$ 
The structural maps satisfy the following properties: \\[1mm] (i)
$r(gh)=r(g)$, $d(gh)=d(h)$ for any pair $(g,h)\in \GR^{(2)}$, and the
partially defined multiplication $\mu$ is associative.\\ (ii)
$d(u(x))=r(u(x))=x$, $\forall x\in \GR^{(0)}$, $u(r(g))g=g$, and
$gu(d(g))=g$, $\forall g\in \GR^{(1)}$, and $u:\Gr0 \to \Gr1$ is
injective.\\ (iii) $r(g^{-1})=d(g)$, $d(g^{-1})=r(g)$,
$gg^{-1}=u(r(g))$, and $g^{-1}g=u(d(g))$.\vspace*{1mm}

By definition, a {\em differentiable groupoid} is a groupoid such
that $\Gr0$ and $\Gr1$ are smooth manifolds with corners, $\Gr0$ 
is smooth, all structural
morphisms are differentiable, and $d$, the domain map, is a
submersion. 

We observe that $\iota$ is a diffeomorphism and hence $d$ is a
submersion if, and only if, $r=d\circ \iota$ is a submersion. Also, it
follows from the definition of submersions of manifolds 
with corners that each fiber $\GR_x=d^{-1}(x) \subset
\Gr1$ is a smooth manifold without corners whose dimension $n$ is constant on each
connected component of $\Gr 0$. The \'etale groupoids considered in
\cite{Brylinski-Nistor1} are extreme examples of differentiable
groupoids (corresponding to $\dim \GR_x =0$). Note that we allow 
$\Gr1$ to be non-Hausdorff. If we want to make more precise the
space of units $\Gr0$ of $\GR$, we say that ``$\GR$ is a differentiable groupoid 
{\em on $\Gr0$}.''

We now recall the definition of a Lie algebroid \cite{Pradines2}.
See also \cite{Rinehart}.

\begin{definition}\label{Lie.Algebroid} A  {\em Lie algebroid}
over a manifold $M$ is a vector bundle $A$ over $M$ together with a
Lie algebra structure on the space $\Gamma(A)$ of smooth sections of
$A$, and a bundle map $q: A \rightarrow TM$, extended to a map between
sections of these bundles, such that:

(i) $q([X,Y])=[q(X),q(Y)]$, and

(ii) $[X, fY] = f[X,Y] + (q(X) f)Y$,

\noindent for all smooth sections $X$ and $Y$ of $A$ and all smooth
function $f$ on $M$.
\end{definition}

The morphism $q$ is called the {\em anchor map}.  Note that we allow
the base $M$ in the definition above to be a manifold with
corners. This is necessary in the most interesting examples related to
singular spaces. We did not include in the definition the condition
that $q(\Gamma(A))$ consist of vector fields tangent to each face of
$M$, but it will be satisfied in all the cases we consider.

To a differential groupoid $\GR$ there is naturally associated a Lie
algebroid $A(\GR)$ as follows \cite{Mackenzie1,Pradines2,Pradines1}.  
Consider
first the vertical tangent bundle of $\GR$ along the fibers of the
domain map $d$:
$$
T_d \GR = \ker (d_*) = \bigcup_{x\in \Gr 0} T \GR_x \subset T\Gr 1.
$$ 
By definition, $A(\GR)$ is the restriction $T_d \GR\vert_{\Gr0}$ of
$T_d \GR$ to the set of units of $\GR$.  The space of $d$--vertical
vector fields invariant with respect to right translations is closed
with respect to the Lie bracket and identifies canonically with
$\Gamma(A)$. Thus we obtain a Lie algebra structure on
$\Gamma(A)$. The action on functions on the base is obtained by
lifting a function on $\Gr0$ to a function on $\Gr1$ via $r$.
Consequently, the anchor map $q:A(\GR) \to T\Gr0$ is obtained by
restricting the differential of $r$ to $A(\GR)$.

If $A$ is a Lie algebroid and $\GR$ is a smooth groupoid such that
$A(\GR) \simeq A$, then we say that $\GR$ {\em integrates} $A$.  Not
every Lie algebroid is integrable (see \cite{al-mo:suites} for an
example).  Nevertheless, it is important to provide examples of
general methods to integrate Lie algebroids.

A differentiable groupoid is called {\em $d$-connected} [respectively,
{\em $d$-simply connected}] if, and only if, each set $\GR_x$ is path
connected [respectively, path connected and simply connected].  As for
Lie groups, all $d$-simply connected differential groupoids with
isomorphic Lie algebroids are isomorphic.

If $\GR$ is a differentiable groupoid,
then there exists a $d$-simply connected groupoid ${\mathcal P}\GR$,
uniquely determined up to isomorphism, with the same Lie algebroid as
$\GR$; it is called the {\em path groupoid} associated to $\GR$. As a
set, ${\mathcal P}\GR$ consists of fixed end-point homotopy classes of
paths $\gamma : [0,1] \to \GR$ such that $d(\gamma(t))$ is constant
and $\gamma(0)$ is a unit. This result is due to Moerdijck in general
(see also \cite{Mackenzie1} for a particular case).

\section{A glueying theorem}

Let $M$ be a manifold with corners. In this paper, by {\em stratified
manifold} we mean a smooth manifold $M$, possibly with corners,
together with a disjoint union decomposition $M=\cup S$ of $M$ by a
locally finite family of submanifolds $S$ without corners, called {\em
open strata}, such that the closure (in $M$) of each stratum $S$ is a
submanifold with corners and each $S$ is contained in a unique open
face of $M$. (This notion is slightly stronger than that of a
stratified space which is also a manifold with corners. Most of our
results however hold in this greater generality.)

Consider a differentiable Lie algebroid $A$ with anchor map $q: A \to
TM$ on a manifold $M$. A stratification $M = \cup S$ of $M$ is called
{\em $A$-invariant} if, and only if, for each point $x \in S$, the
range of $A_x \to T_xM$ (from the fiber of $A$ at the point $x$ to the
tangent space to $M$ at $x$) is contained in $T_x S$ (\ie $q(A_x)
\subset T_x S$). The condition that
\begin{equation}
	q(A_x) \subset T_x S,
\end{equation}
for all $x \in S$, is equivalent to the condition that each local
diffeomorphism of the form $\exp(q(X))$, for some smooth section $X$
of $A$, preserve the strata of $M$, or, to the condition that the
restriction $A_S$ of $A$ to $S$ be a Lie algebroid on $S$, for
each $S$.

Define, for any subset $S\subset M$ of the set of units of $\GR$, the
groupoid
$$
	\GR_S := d^{-1}(S) \cap r^{-1}(S),
$$
called {\em the reduction of $\GR$ to $S$}. If $M = \cup S$ is an
$A$-invariant stratification of $M$ and, moreover, $\GR$ is a
differentiable groupoid on $M$ that integrates $A$, then $d^{-1}(S) =
r^{-1}(S)$ and $\GR_S =d^{-1}(S)$ satisfies $A_S \simeq
A(\GR_S)$. This shows that, in order to integrate $A$, we need to
integrate each restriction $A_S$. The main point of this section is
that this is also enough (Theorem \ref{Theorem.Glueying}).

Recall from \cite{Mackenzie1} that an {\em admissible section} of a
differential groupoid $\GR$ on $M$ is a differentiable map 
$$
	\sigma : M \rightarrow \GR
$$ 
such that $d(\sigma(x))=x$ and the map $M \ni x \to
r(\sigma(x)) \in M$ is a diffeomorphism.  Then $\sigma$ also defines a
diffeomorphism
\begin{equation}\label{eq.l.diff}
	\GR \ni g \rightarrow \sigma g := \sigma(r(g))g \in \GR.
\end{equation}
The main example of an admissible section is 
$$
	\sigma(x)=\exp(X_m) \ldots \exp(X_1)x,
$$ 
for suitable, smooth sections $X_1,\ldots, X_m$ of $A$. We will discuss
this type of admissible sections in more detail below.

A {\em differentiable family} of admissible sections is a family
$\sigma_s: M \to \GR$, $s \in [0,1]$, of maps such that each
$\sigma_s$, $s \in [0,1]$, is an admissible section and the induced
map $[0,1] \times M \ni (s,x) \to \sigma_s(x) \in \GR$ is
differentiable.

\begin{lemma} \label{Lemma.invariant}\
Let $\GR$ be a differentiable groupoid on $M$ and let $\sigma_s: M \to
\GR$ be a differentiable family of local admissible sections. Then
there exists a section $X$ of $A$ such that
\begin{equation}
\label{eq.diff}
	\pa_s f(\sigma_sg)\vert _{s=0}=Xf(\sigma_0g), \;\; g \in \GR.
\end{equation}
\end{lemma}

\begin{proof}
Since the map $[0,1] \ni s \to \sigma_s g:= \sigma_s(r(g))g$ is
differentiable for all $g$, there exists a
vector field $X$ on $\GR$ satisfying \eqref{eq.diff}. We need to check
that $X$ is $d$-vertical and right invariant. We have that
$$
d(\sigma_sg)=d(g),
$$ 
which proves that $X$ is $d$-vertical, by definition. Also,
$$
\sigma_s (gh)=\sigma_s(r(gh))gh=\sigma_s(r(g))gh=(\sigma_s g)h,
$$ 
which proves that $X$ is also right invariant.
\end{proof}

In the following, the section $X \in \Gamma(A(\GR))$, defined in the
above lemma, will be denoted $\pa_s \sigma_s\vert_{s=0}$. We define in
the same way $\pa_s \sigma_s$ for all values of $s$.

We shall repeatedly use the exponential map, and hence we shall to
consider diffeomorphisms obtained by integrating vector fields. Recall
that a smooth vector field $X$ on a manifold $M$ is called {\em
complete} if, and only if, there exists a differentiable map $\phi:
\RR \times M \to M$ such that
$$	
	X(\phi(t,m))=\pa_t(\phi(t,m)),
$$ 
for all $(t,m) \in \RR \times M$. We then define
$$
	\exp(X) m := \phi(1,m).
$$
Note that a complete vector field on a manifold with corners $M$ is
necessarily tangent to each face of $M$.

As for manifolds without corners, it follows from the definition and
basic results on ordinary differential equations that, if $X$ is
complete, then $tX$ is also complete, for all $t \in \RR$, and
$$
	\exp((t+s)X) = \exp(tX)\exp(sX),
$$
for all $s$ and $t$. Consequently, $\exp(X)$ is a diffeomorphism.

\begin{lemma}\ Let $\GR$ be a differentiable
groupoid on $M$ with Lie algebroid $A = A(\GR)$. If $X \in \Gamma(A)$
is a section such that $q(X)$ is complete, then $X$, regarded as a
$d$-vertical vector field on $\GR$, is also complete.
\end{lemma}

\begin{proof} For manifolds without corners,
this is a result from Kumpera and Spencer \cite[Appendix]{KS72}. For manifolds
with corners the proof is the same.
\end{proof}

\begin{proposition}\label{Prop.transport}\
Let $\pi: A \to M$ be a Lie algebroid with anchor map $q: A \to TM$,
and let $X \in \Gamma(A)$ be a section such that $q(X)$ is
complete. Then there exists a uniquely determined isomorphism $E_X :
A \to A$ of Lie algebroids satisfying (i) and (ii):

(i)\ $\pi \circ E_X = \exp(q(X)) \circ \pi$;

(ii)\ $\pa_t E_{tX}(Y)\vert_{t=0} = [X,Y]$.

(iii)\ If, moreover, $\GR$ is a differentiable groupoid integrating $A$ and
$X,Y \in \Gamma(A)$ are both complete, then $E_X$ also satisfies
$$
	\exp(X) \exp(Y) = \exp(E_X(Y)) \exp(X),
$$
as admissible sections.
\end{proposition}

\begin{proof}\ 
If $A$ is integrable, then (i) is Proposition 4.1.(v), (ii) is
Proposition Proposition 4.11.(iii), and (iii) is Proposition 4.11.(ii)
of \cite{Mackenzie1}. The proof of (i) and (ii) in the general case
are the same.
\end{proof}

A family $Y_1, \ldots, Y_n$ of smooth sections of $A$ is a
{\em local basis} of $A$ at $y \in M$ if $Y_1(y), \ldots, Y_n(y)$ is
a basis of $A_y$. If $t =(t_1,\ldots, t_n) \in \RR^n$ and $Y_1,
\ldots, Y_n$ are sections of $A$, then we denote 
$$
	\Exp(t, Y):=\exp(t_1Y_1)\exp(t_2Y_2) \ldots  \exp(t_nY_n).
$$  
Also, let 
$$
	B_{\epsilon} = \{ t =(t_1,\ldots, t_n) \in \RR^n,\; 
	\sum_{i=1}^n |t_i| <
	\epsilon\}.
$$

Recall that a differentiable local groupoid ${\mathcal U}$ on $M$
satisfies all the axioms of a differential
groupoid, except that the multiplication $gg'$ is not defined for all
pairs $(g,g')$ such that $d(g) = r(g')$, but only in a neighborhood
${\mathcal U}_2$ of the diagonal in the set $\{(g,g'), d(g) =
r(g')\}$. See \cite{vanEst,Pradines1} for details.

\begin{lemma} \label{Lemma.Coord.Sys}\
Let $\GR$ be a local differential groupoid on $M$, and let $\sigma: M
\to \GR$ be an admissible section with $r(\sigma(x_0))=y_0$. Also, let
$Y_1,\ldots, Y_n$ be a local base of $A$ at $y_0$. Then, for some
small $\epsilon > 0$ and a relatively compact open neighborhood $U$ of
$x_0$ in $M$, the map
\begin{equation}
\label{eq.coord.syst}
	\psi_Y^\sigma : \RR^n \times M \ni (t,x) \to \Exp(t, Y)
	\sigma(x) \in \GR
\end{equation}
is a diffeomorphism from $B_{\epsilon} \times U$ to a neighborhood of
$\sigma(x_0)$ in $\GR$. 

If $\sigma=\exp(X_m) \ldots \exp(X_1)$, for
some integrable $X_1,\ldots,X_m \in \Gamma(A)$, then we denote
$\psi_Y^\sigma=\psi_Y^X$.
\end{lemma}

\begin{proof}\ If $\GR$ is a differentiable groupoid (not just a local one), 
this is Proposition 4.12 in \cite{Mackenzie1}. For local groupoids the
proof is the same.
\end{proof}

Let $\GR$ be a differentiable groupoid, and let $X_1,\ldots, X_m$ be
smooth complete sections of $A(\GR)$. Also, let $Y_i$ be smooth,
complete sections of $A(\GR)$ that form a local basis at some point
$x_0 \in M$, as above. For simplicity, we shall sometimes assume that
the $Y_i$ are compactly supported, which is a stronger assumption than
completeness. Since the admissible section 
$$
	\sigma(x)=\exp(X_m) \ldots \exp(X_1)x
$$ 
defines a diffeomorphism $\sigma:\GR \to \GR$, equation
\eqref{eq.l.diff}, we obtain that the map
\begin{equation}
\label{eq.coord.systXY}
	\phi_Y^X : \RR^n \times M \ni (t,x) \to \exp(X_m)
	\ldots \exp(X_1) \Exp(t, Y) x \in \GR
\end{equation}
is also a diffeomorphism from a set of the form $B_{\epsilon} \times U$
to its image, for some small $\epsilon >0 $
and some small open subset $U \subset M$. The maps $\phi_Y^X$ are slightly
more convenient to work with, in what follows, than the maps $\psi$ of the
previous lemma.

{\em Throughout the rest of this section, $M$ will be a smooth
manifold and $A \to M$ will be a Lie algebroid. Moreover, $M = \cup S$ is an
$A$--invariant stratification of $M$, $\GR_S$ is a differentiable
groupoid integrating $A_S$, and $\GR = \cup \GR_S$ (disjoint union).}

On $\GR=\cup \GR_S$ we can then define uniquely a natural groupoid
structure on $\GR$ such that the structural morphisms of each $\GR_S$
are obtained from those of $\GR$ by restriction.  In particular, if
two arrows $g \in \GR_S$ and $g' \in \GR_{S'}$ are composable, then $S
= S'$, also, each $\GR_S$ is a subgroupoid of $\GR$.

A crucial observation is that we need not use the full differentiable
structure on $\GR$ to define the maps $\psi_Y^X$ of Lemma
\ref{Lemma.Coord.Sys}, or the maps $\phi^X_Y$ of equation
\eqref{eq.coord.systXY}, for that matter. To define $\phi_Y^X$, it is
enough to use only the smooth structure on each $\GR_S$. This will
allow to extend the definition of $\phi_Y^X$ to our case, when the
groupoid $\GR$ is obtained by glueying differential groupoids. This can
be done as follows.

Let $M = \cup S$ be an $A$-invariant stratification of $M$. Also, let
for each $S$ $\GR_S$ be a differentiable groupoid integrating
$A_S=A\vert_S$ and $\GR = \cup \GR_S$, as above.  If $X_i$ and $Y_j$
are sections of $A$ on $M=\Gr0$ such that the vector fields $q(X_i)$
and $q(Y_j)$ are integrable, then the restriction of these vector
fields to each strata is again integrable, and hence the restriction
of $X_i$ and $Y_j$ to each strata are integrable as vertical vector
fields on $\GR_S$. Then the map $\phi_Y^X$ is defined on each $\RR^n
\times S$ (with values in $\GR_S$), by glueying these maps, we obtain
the desired definition of $\phi_Y^X : \RR^n \times M \to \GR$, for
$\GR = \cup \GR_S$.

Since the maps $\phi_Y^X$ and $\psi_Y^X$ play an important role in what 
follows, we now spell out their properties in more detail.

\begin{lemma}\label{Lemma.new1}\
(i) The maps $\phi_Y^X$ and $\psi_Y^X$ are related by 
$$
	\phi^X_Y=\psi^X_{Y'},
$$
where $Y'_j=E_{X_m}\circ E_{X_{m-1}} \circ \ldots \circ E_{X_1}(Y_j)$.

(ii) Let $Y_j''=-Y_{n+1-j}$ and $X_i''=-X_{m+1-i}$. Then 
$$
	\phi^X_Y(x)^{-1} = \psi^{X''}_{Y''}(\alpha(x)),
$$
where $\alpha$ is the diffeomorphism $\exp(q(X_m))\ldots \exp(q(X_1))$.
\end{lemma}

\begin{proof}\  
Because each $\GR_S$ is a differential groupoid, (i) follows from
Proposition \ref{Prop.transport} on each $\GR_S$. From this we obtain
the desired relation everywhere on $\GR$.

(ii) follows from the relation
$$
	(\exp(X)x)^{-1}=\exp(-X)( \exp(q(X))x ),
$$
valid on each $\GR_S$, and hence everywhere on $\GR$. 
\end{proof}

We let as above $X_1, \ldots, X_m,Y_1,\ldots,Y_n$ be integrable sections
of $A$. We shall also use the following lemma.

\begin{lemma}\label{Lemma.new2}\ 
(i) Fix $t \in \RR^n$ and $x_0 \in M$. Then we can find $\epsilon >0$, a
neighborhood $U$ of $x$ in $M$, and a differentiable map
$$
	\tau: B_{\epsilon} \times U \to \RR^n \times M, \;\tau(0,x_0)=(0,x_0),
$$
which is a  diffeomorphism onto its image, such that
\begin{equation*}
	\phi_Y^X(t+s,x)=\phi_{Y}^{X'}(\tau(s,x)),
\end{equation*}
on $B_\epsilon \times U$, 
where $X'_j=X_j$, for $j=1,\ldots,m$, and $X'_{j+m}= t_jY_j$, for $j=1,\ldots,n$.

(ii) For all $Y_1,\ldots Y_n $, $Y'_1, \ldots,Y'_n$, and $x_0 \in M$,
there exist $\delta > \epsilon >0$, an open neighborhood $U$ of $x_0$,
and a differentiable map
$$
	m_1 : B_{\epsilon} \times B_\epsilon \times U \to B_{\delta} \times M,\;
	m_1(0,0,x)=(0,x),$$ 
such that
$$
	\Exp(t,Y)x (\Exp(t',Y')x)^{-1} = \Exp(m(t,t',x)).
$$
\end{lemma}

\begin{proof}\ We begin by writing $\exp((t_i + s_i) X_i)=\exp(t_iX_i)
\exp(s_iX_i)$. Then, using Proposition \ref{Prop.transport}, we can find
sections $Y_i^s$ of $A$, $s \in \RR$, $Y_i^0=Y_i$, depending
smoothly on $s$, such that
$$
	\phi_Y^X(t+s,x)=\phi_{Y}^{X'}(\Exp(s,Y^s)x).
$$
To complete (i), we now use a local groupoid $\mathcal U$ integrating
$A$. The existence of such a $\mathcal U$ is ensured by
\cite{Pradines1}. By replacing $\mathcal U$ with an open neighborhood 
of $M$, if necessary, we may identify $\mathcal U$ with a subset of $\GR$, 
using the exponential map. Then, for small $\epsilon$ and a relatively
compact neighborhood $U$ of $x_0$, the maps 
$$
	B_\epsilon \times U \ni (x,s) \to f_1(x,s):=\Exp(s,Y^s)x
$$ 
and 
$$
	B_\epsilon \times U \ni (x,s) \to f_2(x,s):=\Exp(s,Y)x
$$ 
are diffeomorphisms onto neighborhoods of $x_0$ in $\mathcal U$, by
Lemma \ref{Lemma.Coord.Sys}.  The desired map $\tau$ is obtained from
$f_2^{-1} \circ f_1$.

The proof of (ii) is similar, using the differentiability of
the multiplication in a local groupoid.
\end{proof}

The above result suggest to introduce the following family of
maps.\vspace*{1mm}

{\sc The family} $\Phi$:\ \ Let $f: V_0 \to V$ be a diffeomorphism
(\ie coordinate chart) from an open subset of $\RR^l$ to an open,
relatively compact subset of $\Gr0=M$. Because the partition $M =\cup
S$ is locally finite, we can find, using Lemma \ref{Lemma.Coord.Sys},
an $\epsilon > 0$ and sections $X_i,Y_j \in \Gamma(A)$, for which
$\phi_Y^X$ defines a diffeomorphism from $B_{\epsilon} \times (V \cap
S) \to \GR_S$, for all $S$ such that the intersection $S \cap V$ is
not empty. The family $\Phi$ then consists of all maps of the form
\begin{equation}\label{eq.Phi}
	\varphi(t,y)=\phi_Y^X(t,f(y)) : B_\epsilon \times V_0 \to \GR.
\end{equation}

Recall that a {\em differentiable atlas} on a set $M_0$ is a family
of injective maps 
$$
	\varphi : V_\varphi \to M_0,
$$
defined on an open subset of $\RR^l$, for
some fixed $l$, such that $\varphi(V_\varphi)$ is a covering of $M_0$,
$\varphi^{-1}(\varphi_1(V_{\varphi_1})$ is an open subset
$V$ of $V_{\varphi}$, and the map $\varphi_1^{-1}\varphi$ is differentiable
on $V$.

We are ready now to prove the following theorem.

\begin{theorem} \label{Theorem.atlas}\
Let $M = \cup S$ be an $A$-invariant stratification of $M$ and, for
each $S$, let $\GR_S$ be a $d$-connected differentiable groupoid on
$M$ integrating $A_S=A\vert_S$. Then $\GR$ has a differentiable
structure making it a differentiable groupoid with $A(\GR)\simeq A$
if, and only if, the family $\Phi$, consisting of the maps in equation
\eqref{eq.Phi}, is a differentiable atlas.
\end{theorem}

\begin{proof}\  Suppose first that the family $\Phi$ is a differentiable
atlas. We begin by showing that the groupoid structure on $\GR$,
induced from the groupoids $\GR_S$, is compatible with the
differentiable structure defined by $\Phi$.  That is, we need to check
that the structural morphisms are differentiable. We now check this.

First, the domain map is differentiable and submersive because
$d(\phi_Y^X(t,x))=x$, for all $X$, $Y$, and $x \in M$. Next, it is
enough to check that the map $(g',g)\to g'g^{-1}$, defined on 
$$
	\{(g',g), d(g')=d(g)\},
$$ 
is differentiable. Let $(g',g)$ be such that $d(g)=d(g')=x_0$, and let
$\phi_Y^X : \RR^n \times M \to \GR$ and $\phi_{Y'}^{X'} : \RR^n \times
M \to \GR$ be two maps such that $\phi_Y^X(t,x_0)=g$ and
$\phi_{Y'}^{X'}(t',x_0)=g'$, and such that their restrictions to some
small set of the form $B_\epsilon \times U$, $t,t' \in B_\epsilon$, is
in the family $\Phi$. We need to show that the induced map
\begin{equation}
	\mu_1:B_\epsilon \times B_\epsilon \times U \ni (t,t',x) 
	\rightarrow \phi_{Y'}^{X'}(t',x)(\phi_Y^X(t,x))^{-1} \in \GR,
\end{equation}
is differentiable. It is enough to prove this in a small neighborhood
of $x_0$. Because $\Phi$ forms an atlas, using Lemma
\ref{Lemma.new2}(i), we see that we may assume $t=t'=0$, eventually by
changing $X,X',Y,$ or $Y'$. The differentiability of $\mu_1$ then
follows by combining Lemma \ref{Lemma.new2}(ii), and Lemma
\ref{Lemma.new1}(i). This is enough to conclude that $\GR$ is a
differentiable groupoid whenever $\Phi$ is an atlas.

Let $A(\GR)$ be the Lie algebroid of $\GR$ corresponding to the
differentiable structure defined by $\Phi$. We need to check that
$A(\GR) \simeq A$. We have that $A(\GR)\vert_{S} \simeq A_S$, by
construction, and hence we can identify as a set $A(\GR)$ with the
disjoint union of the restrictions $A_S$. The two differentiable
structures on $\cup A_S$ (the first induced from $A$ and the second
induced from $A(\GR)$) are the same because the differential of the
map $\phi_Y^\sigma$ of Lemma \ref{Lemma.Coord.Sys}, $\sigma = id$,
canonically identifies $A\vert_U$ and $A(\GR)\vert_U$, if $U$ is as in
that lemma. The Lie algebra structures on $\Gamma(A)$ and
$\Gamma(A(\GR))$ also coincide because the two possible brackets of
two vector fields coincide on each strata $S$, and hence they coincide
everywhere. We have thus proved that, if $\Phi$ forms an atlas, then
$\GR$ is a differentiable groupoid with $A(\GR) \simeq A$.

Conversely, suppose now that $\GR$ is endowed with a differentiable 
structure, and let 
$g \in \GR_S$. Since $\GR_S$ is $d$-connected, we can choose
vector fields $X_1, \ldots, X_m$ such that
$$
	g = \exp(X_m) \ldots \exp(X_1) x_0,
$$
for some $x_0 \in S$. If $Y_1, \ldots, Y_n$ are chosen to form a basis
of $A_{x_0}$, as in Lemma \ref{Lemma.Coord.Sys}, then the map
$\phi_Y^X$ must be a diffeomorphism of a set of the form $B_\epsilon
\times U$ onto an open neighborhood of $g$, for some open subset $U
\subset M$.  We obtain that if $\GR = \cup \GR_S$ is a differential
groupoid such that $A(\GR) \simeq A$, then the family $\Phi$ is an
atlas.
\end{proof}

In particular, we immediately obtain from the above theorem and its
proof that the differentiable structure on $M$ and the groupoid
structure on $\GR$ uniquely determine the differentiable structure on $\GR$
satisfying $A(\GR) \simeq A$. Indeed, the differentiable structure on 
$\GR$ is determined by the family $\Phi$.

We now turn to the main theorem of this paper, Theorem
\ref{Theorem.Glueying}, which shows that the assumptions of Theorem
\ref{Theorem.atlas} are satisfied provided that the groupoids $\GR_S$
are $d$-simply connected. We shall need an extension of the concept of
differentiable family of sections to our case, $\GR = \cup \GR_S$,
when $\GR_S$ are smooth groupoids but $\GR$ is not endowed with any
differentiable structure. 

A {\em differentiable family} of admissible
sections is a family $\sigma_s: M \to \GR$, $s \in [0,1]$, of maps
such that each $\sigma_s$, $s \in [0,1]$, is an admissible section,
the induced map $[0,1] \times S \ni (s,x) \to \sigma_s(x) \in \GR_S$
is differentiable for each $S$, and the sections
$\pa_s(\sigma_s\vert_S)$ of $A_S$ can be glued to a smooth section of
$A$ on $[0,1] \times M$. This definition of a differentiable family of
smooth section is thus very similar to the corresponding definition in
the case when $\GR$ has a smooth structure, except that we replace the
condition on the smoothness of the map $[0,1] \times M \to \GR$ by the
existence and smoothness of the derivatives $\pa_s \sigma_s$ on 
$[0,1] \times M$.

The following lemma is an important technical part of the proof of
Theorem \ref{Theorem.Glueying}. It achieves a continuous and smooth
deformation of a local admissible section of $\GR$ to a local identity
section. We denote by $\Gamma_c(E)$ the space of compactly supported,
smooth sections of a vector bundle $E$.

\begin{lemma} \label{Lemma.simply.connected}\
Let $\GR = \cup \GR_S$ be a union of $d$-simply connected
differentiable groupoids with $A(\GR_S)=A_S$, as in the statement of
Theorem \ref{Theorem.atlas}, and let $x_0 \in M$.  Suppose that
$X_1, \ldots, X_m \in \Gamma_c(A)$ satisfy
$$
	\exp(X_m) \ldots \exp(X_1)x_0=x_0.
$$ 
Then we can find
a differentiable family of admissible sections
$\sigma_s: M \to \GR$ such that:

(i) \ $\sigma_1 (x) = \exp(X_m) \ldots \exp(X_1)x$, on $M$;

(ii) \ $\sigma_0(x)=x$, for all $x \in M$; \\ and, most importantly,

(iii) \ $\sigma_s(x_0)=x_0,$ for all $s \in [0,1]$.
\end{lemma}

This lemma remains true if we replace the condition that $X_i$
be compactly supported by the condition that they be integrable.

\begin{proof}
Consider the curve 
$$
	\phi(t)=\exp(t'X_k)\exp(X_{k-1}) \ldots \exp(X_1)x_0,
$$
where $t'=(mt-k +1)$ and $k$ is chosen such that $0 \le t' \le 1$. By
assumption, $\phi$ is a closed curve on $\GR_{x_0}$, and hence,
by the assumption that $\GR$ is $d$--simply connected, we can
continuously deform this curve to the constant curve $x_0$, within
$\GR_{x_0}$, through closed curves based at $x_0$. More precisely, we
can find 
$$
	\eta:[0,1] \times [0,1] \to \GR_{x_0}
$$ 
such that $\eta(t,1)=\phi(t)$ and $\eta(t,0)=\eta(0,s)
=\eta(1,s)=x_0$, for all $t,s \in [0,1]$.

By an approximation argument, we can assume that $\eta(\frac{k}{m},s)$
depends smoothly on $s$, for each integer $k$.  Moreover, after
replacing $m$ by a large multiple $lm$ and each $X_k$, $k = 1, \ldots, m$,
by $l^{-1}X_k$ repeated $l$ times, we can assume that there exist
compactly supported
sections $X_k^s \in \Gamma(A)$, depending smoothly on $s$, such
that $X_k^1=X_k$, $X_k^0=0$, and
$$
	\eta(\frac{k}{m},s)=\exp(X_k^s)\eta(\frac{k-1}{m},s).
$$
The desired deformation is obtained by letting
$$
	\sigma_s(x)=\exp(X_m^s) \ldots \exp(X_1^s)x.
$$
To see that $\sigma_s$ is a smooth family of admissible sections, we
use Proposition \ref{Prop.transport} for each $\GR_S$ and the fact
that each $\exp(X_k^s)$ is a smooth admissible section.
\end{proof}

We continue to assume that $\GR =\cup \GR_S$ is as in the statement
of Theorem \ref{Theorem.Glueying}.

\begin{proposition} \label{Prop.diffeo}\
Let $\sigma_s$, $s \in[0,1]$, be a differentiable family of local
admissible sections. Assume that $\sigma_s(x_0)=x_0$, for all $s \in [0,1]$, 
and that $\sigma_0(x)=x$, for all $x \in M$. Let
$Y_1, \ldots, Y_n$ be a local basis of $A$ at $x_0$. Then there exist
$\epsilon > \delta > 0$, a neighborhood $U$ of $x_0$ in $M$, and a
differentiable map
$$
	\tau: [0,1] \times B_{\delta} \times U \to B_\epsilon
$$ 
such that:\ (i)\  $\tau(s, 0, x_0)=0$,\ (ii)\ 
for each fixed $s\in [0,1]$ and $x\in
U$, the map $B_{\delta} \ni t \to \tau(s,t,x) \in B_{\epsilon}$ 
is a diffeomorphism onto its image, and
\begin{equation*} 
	(iii) \qquad \sigma_s \Exp(t, Y) x = \Exp(\tau(s,t,x),Y) x,
\end{equation*}
for all $s, t \in [0,1]$, and $x\in U$.
\end{proposition}

\begin{proof}\ Let ${\mathcal U}$ be a local groupoid integrating
$A$. Choose $\epsilon > 0$ and a neighborhood ${U}_1$ of $x_0$ small
enough such that
$$
	\eta: B_{2\epsilon} \times U_1 \to {\mathcal U}, \;\;
	\eta(t,x) =\Exp(t, Y)x,
$$
is a diffeomorphism onto an open neighborhood $V$ of $x_0$ in
${\mathcal U}$.  Let $U_0 \subset U_1$ be a compact neighborhood of
$x_0$. Because $K:=\exp(\overline{B}_{\epsilon} \times U_0)$ is a
compact subset of the open set $V$ and $\sigma$ is continuous,
$\sigma_s(x_0)=x_0$, we can find a neighborhood $U$ of $x_0$ such that
$\sigma_s k$ is defined in ${\mathcal U}$ whenever $k \in K$, and
$d(k) \in U$, and such that
$$
	\sigma_s(U_1) K \subset V.
$$
Then we can define 
$$
	\tau(s,t,x)=\eta^{-1}(\sigma_s\Exp(t, Y)x).
$$
Because $\sigma_s\Exp(t, Y)x \in V \subset {\mathcal U}$ 
and $\eta$ is a diffeomorphism, we obtain that $\tau$ is smooth also.
\end{proof}

\begin{theorem} \label{Theorem.Glueying}\
Let $A$ be a Lie algebroid on a manifold with corners $M$.  Suppose
that $M$ has an $A$-invariant stratification $M = \cup S$ such that,
for each stratum $S$, the restriction $A_S$ is integrable,
then $A$ is integrable.

More precisely, let $\GR_S$ be $d$-simply connected differential
groupoids such that $A(\GR_S) \simeq A_S$.  Then the disjoint union
$\GR = \cup \GR_S$ is naturally a differentiable groupoid such that
$A(\GR) \simeq A$.
\end{theorem}

\begin{proof}\  By Theorem \ref{Theorem.atlas}, we see that 
we it is enough to show that the family $\Phi$, defined using the maps
$\phi_Y^X$ of equation \eqref{eq.coord.systXY}, is an atlas. Let
$\phi_Y^X$ and $\phi_{Y'}^{X'}$ be two such maps, defined on
$B_\epsilon \times U$ and, respectively, on $B_{\epsilon'} \times U'$,
such that their images intersect. Let $g \in \GR$ be an element of
this intersection. As in the proof of Theorem \ref{Theorem.atlas}, we
can arrange that $\phi^X_Y(0,x_0)=\phi_{Y'}^{X'}(0,x_0)=g$. Then
$$
	\phi^X_Y(0,x_0)=\exp(X_m) \ldots \exp(X_1) x_0=
	\exp(X'_{m'}) \ldots \exp(X'_1) x_0=\phi_{Y'}^{X'}(0,x_0),
$$
and hence 
$$
	x_0= \exp(-X'_1) \ldots \exp(-X'_{m'}) \exp(X_{m}) \ldots \exp(X_1)
	x_0,
$$
so we may also assume that $X'_i=0$ (and hence $g=x_0$). By replacing
$V$ and $V'$ with some smaller, relatively compact neighborhoods, if
necessary, we may further assume that the sections $X_i$ are compactly
supported. Let
$$
	\sigma(x) = \exp(X_{m}) \ldots \exp(X_1)x.
$$  
Since each $\GR_S$ is $d$-simply connected, by Lemma
\ref{Lemma.simply.connected} we can find a smooth family $\sigma_s:M
\to \GR$, $s \in [0,1]$, such that $\sigma_0$ is the identity,
$\sigma_1=\sigma$, the restriction to each groupoid $\GR_S$ is
differentiable. The smoothness of the family $\sigma_s$ in this case
is reduces to showing that ${\pa_s}\sigma_s$ is a smooth section of
$A$ over $M$. Let $\tau$ be the map defined in Proposition
\ref{Prop.diffeo} using the admissible sections $\sigma_s$, and let
$\tau_1(t,x) = \tau(1,t,x)$.  Then
$$
	\phi^X_Y(t,x)=\phi^{X'}_{Y'}(\tau_1(t,x),x),
$$
in a small neighborhood of $x_0$.
Since $\tau_1$ is a local diffeomorphism for each fixed $x$, 
this proves the theorem.
\end{proof}

\section{Applications to foliations}

Recall that a Lie algebroid $\pi : A \to M$, with anchor map $q : A
\to TM$, is called {\em regular} if, and only if, the range of $q$ has
locally constant rank. Then the sections of $q(A)$ are the sections of
the tangent bundle to a foliation ${\mathcal F}_A$ whose tangent
bundle, a sub-bundle of $TM$, is denoted $T{\mathcal F}_A$. If,
moreover, there exists a morphism $\rho:\Gamma(T{\mathcal F_A}) \to
\Gamma(A)$ such that $q \circ \rho = id$, then we may assume that $A$
is a {\em semi-direct product}.

Also, if $A$ is regular, the kernel of $q$ is a bundle of Lie
algebras. See also \cite{Kubarski,Weinstein1}. Recall that a result of
Douady and Lazard from \cite{Douady-Lazard} states that every bundle
of Lie algebras is integrable.  In particular the kernel $\ker (q)$ is
integrable. (I am greatful to Alan Weinstein for pointing out this
reference to me.) 

We say that $A$ is (isomorphic to) the semi-direct product $\ker(q)
\rtimes q(A)$ if there exists a morphism of Lie algebras $\rho:
\Gamma(q(A)) \to \Gamma(A)$ that is a right inverse to $q$.

\begin{proposition} \label{Prop.integrable}
Let $A$ be a regular algebroid with anchor map $q:A \to TM$.  Assume
that $A$ is the semi-direct product $\ker(q) \rtimes q(A)$.  Then $A$
is integrable.
\end{proposition}

Before proceeding to the proof, we first introduce some terminology and
make some comments.
If $E \to M$ is a vector bundle on a foliated manifold $M$, then a
leafwise connection on $E$ is a linear map 
$$
	\nabla: \Gamma(E) \to \Gamma(E \otimes T^*{\mathcal F}),
$$
satisfying the Leibnitz identity.  Here ${\mathcal F}$ is the
foliation of $M$, $T{\mathcal F}$ is the tangent bundle to this
foliation, and $T^*{\mathcal F}$ is the dual of this bundle.  The
equivalent definitions of a connection in terms of parallel transport
or equivariant splittings of tangent spaces of principal bundles
extend to the ``leafwise'' setting also.

The structure of regular, semi-direct product Lie algebroids is as 
follows.  First, there exists a {\em leafwise flat} connection $\nabla$,
$\nabla_Z(X)=[\rho(Z),X],$ on $\ker (q)$, which preserves its Lie
bundle structure, that is
$$
	\nabla_Z \nabla_{Z'} -\nabla_{Z'} \nabla_{Z} = \nabla_{[Z,Z']}, \quad 
	\Mand
$$
$$
	\nabla_Z([X,Y]) = [\nabla_Z(X),Y] + [X,\nabla_Z(Y)],
$$
for all $Z, Z' \in \Gamma(T{\mathcal F}_A)$ and $X,Y \in \Gamma(\ker
(q))$. Then
$$
A \simeq T{\mathcal F}_A \oplus \ker (q)
$$ 
as vector bundles, with anchor map given by the projection onto the
first component, and with Lie bracket on $\Gamma(A)$ defined by
$$
[(Z,X),(Z',X')]=\big([Z,Z'], \nabla_Z(X') - \nabla_{Z'}(X) + [X,X']\big),
$$
for all $Z,Z' \in \Gamma(T{\mathcal F}_A)$ and $X,X' \in \Gamma(\ker
(q))$. 
 
With these comments, we are now ready to prove Proposition
\ref{Prop.integrable}.

\begin{proof} 
Recall that $\mathcal{PF}$, the {\em path groupoid} of the foliation
${\mathcal F}$, consists of fixed end point homotopy classes of paths
$\gamma$ that are fully contained in a single leaf, with respect to
homotopies {\em within} that leaf. As explained above, the morphism
$\rho$ defines a leafwise flat connection on $\ker(q)$ that preserves
the Lie bracket. Let ${\mathcal K}$ be a $d$-simply connected groupoid that
integrates $\ker (q)$.

We define then a groupoid $\GR$ that integrates $A$ as follows. As a
smooth manifold,
$$
\GR =\{ (g,\gamma) \in {\mathcal K} \times \mathcal{PF},
d(g) = \gamma(1) \}.
$$
To define the multiplication, observe first that the leafwise flat 
connection on $\ker(q)$ defines a parallel transport map
$$
\rho(\gamma) : \ker(q)_{\gamma(0)} \to \ker(q)_{\gamma(1)},
$$
which is a Lie algebra isomorphism for any path $\gamma$ fully
contained in a leaf. Since ${\mathcal K}_x$ is simply connected for
each $x$, we obtain by exponentiation a group morphism
$$
\rho(\gamma):{\mathcal K}_{\gamma(0)}
\to {\mathcal K}_{\gamma(1)}.
$$
That is, the leafwise flat connection on $\ker(q)$ lifts to a leafwise
flat connection 
on ${\mathcal K}$ that preserves the Lie group structure on the 
fibers.

We are ready now to define the groupoid structure 
 on $\GR$.  Note first that $d(g)=r(g)$ for $g \in
{\mathcal K}$. Then $d(g,\gamma)=\gamma(0)$, $r(g,\gamma)=r(g)$, and
the product on $\GR$ is given by the formula
$$
(g,\gamma)(g',\gamma')=(g \rho(\gamma)(g'),\gamma \gamma'),
$$
where the composition of paths is given by concatenation. The flatness
of $\nabla$ gives that $\rho(\gamma\gamma')=\rho(\gamma)\rho(\gamma')$,
which guarantees the associativity of the product.
\end{proof}

In certain cases we get integrability without assuming that $A$ is
a semi-direct product.

\begin{proposition}\ \label{Prop.integrable2}
Let $A$ be a regular Lie algebroid with anchor map $q:A \to TM$ such
that $\ker (q)$ is a bundle of {\em semisimple} Lie algebras. Then $A$
is integrable.
\end{proposition}

\begin{proof}
We may assume that $M$ is connected. Since semisimple Lie algebras are
rigid, all fibers of $\ker(q)$ will be isomorphic Lie algebras.
Fix one of these algebras and denote it by ${\mathfrak g}$. Also, let
$G_0$ be the group of automorphisms of ${\mathfrak g}$. Then, if we
define 
$$
P = \cup_x \operatorname{Iso}({\mathfrak g}, \ker(q)_x),
$$
(the fibers are the sets of Lie algebra isomorphisms ${\mathfrak g}\to
\ker(q)_x$), we obtain a $G_0$-principal bundle on $M$, which acquires
by pull-back a foliation ${\mathcal F}$ of the same codimension as
${\mathcal F}_A$.  The path groupoid $\mathcal{PF}$ of the foliation
${\mathcal F}$ has an induced free action of $G_0$, and we define the
groupoid $\GR$ by $\GR = \mathcal{PF}/G_0$.  The composition of two
paths in $\GR$ is obtained by choosing composable liftings in
$\mathcal{PF}$. Because the Lie algebra of $G_0$ is $\mathfrak g$, we
obtain that $\GR$ integrates $A$.
\end{proof}

With the above results in mind, we now define $A$-invariant regular
stratifications.

\begin{definition}\
Let $\pi : A \to M$ be a Lie algebroid with anchor map 
$q : A \to TM$ on the manifold with corners $M$. An invariant  stratification 
$M = \cup S$ is called {\em regular} if, and only if,
the restriction $A_S:=A\vert_S$ is regular for each $S$.
\end{definition}

Although almost all interesting Lie algebroids have invariant regular
stratifications, this is not true in general.  Consider, for example,
a closed subset $B \subset \RR^n$ with empty interior, which is not a
manifold. Let $\phi \ge 0$ be a smooth function that vanishes exactly
on $B$, and let ${\mathcal V}(\RR^n)$ be the Lie algebra of vector
fields on $\RR^n$.  Since $\phi{\mathcal V}(\RR^n)$ is a free
$\CI(\RR^n)$ module, using the Serre-Swan theorem, we can define a
vector bundle $A$ such that $\Gamma(A) = \phi {\mathcal
V}(\RR^n)$. Moreover,
$$
[\phi X, \phi Y]=\phi\big( X(\phi) Y - Y(\phi) X + \phi [X,Y]\big),
$$
for all vector fields $X$ and $Y$ on $\RR^n$, which shows
that $A$ is a Lie algebroid. It is not difficult to see that $A$ has
no invariant regular stratification.

Although not all Lie algebroids have regular stratifications,
this notion is useful because it is easier to
integrate regular algebroids than general Lie algebroids, 
and we know that in order to integrate a
Lie algebroid, it is enough to integrate it over each strata
(Theorem \ref{Theorem.Glueying}).

\begin{theorem}\ \label{Theorem.application}
Let $A$ be a Lie algebroid with anchor map $q : A \to TM$ on a
manifold $M$ with a regular $A$-invariant stratification $M =\cup S$.
Assume, for each $S$, that either $A_S$ is the semi-direct product
$\ker(q)_S \rtimes q(A_S)$, or that $\ker (q)\vert_S$ is a bundle of
semisimple Lie algebras. Then $A$ is integrable.
\end{theorem}

\begin{proof}\ Using Proposition \ref{Prop.integrable} or Proposition
\ref{Prop.integrable2}, we see that each of the algebroids $A_S$,
obtained by restricting $A$ to the stratum $S$, is integrable.  By
Theorem \ref{Theorem.Glueying}, the Lie algebroid $A$ is then
integrable.
\end{proof}

Let us consider now in greater detail a class of examples that is
useful for the construction of pseudodifferential operators on complex
algebraic varieties. It is possible to associate several ``natural''
algebroids to a complex algebraic variety, so the following
constructions lead to a family of algebras of pseudodifferential
operators associated to a complex algebraic variety. The details of
this construction will be presented in a future paper. We do not know
at this point which of the many algebras that one obtains is the
``right'' algebra to associate to a complex algebraic manifold, but we
hope to address this question in a future paper.

Let $M$ be a manifold with corners, such that
each hyperface (\ie face of maximal dimension) $H \subset M$ is given
by $H = \{ x_H = 0 \}$ for some function $x_H \ge 0$ with $dx_H \not =
0 $ on $H$, (\ie $x_H$ is a defining function of $H$.)  If $F\subset
M$ is an arbitrary face of $M$ of codimension $k$, then $F$ is an open
component of the intersection of the hyperfaces containing it.  The
set $x_1,\ldots,x_k$ of defining functions of these hypersurfaces are
called the defining functions of $F$; thus $F$ is a connected
component of $\{x_1=x_2 =\ldots =x_k=0\}$.

We first introduce the class of Lie algebroids we are interested in on
a manifold with corners $M$. We call these algebroids
``quasi-homogeneous,'' and we now proceed to construct them. First we
need some related definitions that will make it easier to describe our
settings.

\begin{definition}\
A {\em Lie flag} on a manifold $M$ is an increasing finite sequence of
sub-bundles 
$$
	E_0 \subset E_1 \subset \ldots \subset E_l \subset E_\infty := TM
$$
such that $ [\Gamma(E_i) , \Gamma(E_j)] \subset \Gamma(E_{i+j})$.
\end{definition}

It follows from the definition that, if $E_0\subset \ldots \subset E_l
\subset TM$ is a Lie flag on $M$, then the bundles $E_0$ and $E_l$
are integrable. We do not assume the above inclusions to be strict.

We now describe the type of behavior we want for the vector fields
that are sections of a quasi-homogeneous Lie algebroid. Let $p_i$
denote the projection onto the $i$th component of a product.

\begin{definition}\
Let $H = \{x_H=0\} \subset M$ be a hyperface, and let 
$$
\phi_H=(\pi_H,x_H): V_H \to H \times [0,\epsilon)
$$ 
be a diffeomorphism defined in a neighborhood $V_H$ of $H$ in $M$.
Also, let 
$$
	E_0^H \subset \ldots \subset E_{l_H}^H \subset TH
$$ 
be a Lie flag on $H$ and $d_H \in \NN \cup \{\infty\}$.  Then a vector
field $X$ on $V_H$ is called $(E_i,\phi_H,d_H)$--{\em adapted} if, and
only if,
\begin{equation}\label{eq.adapted}
 	X \in \sum_{j=0}^{l_H} x_H^j \Gamma(\pi_H^*(E_j^H)\vert_{V_H}) +
 	\CI(V_H) x_H^{d_H +1} \pa_{x_H},
\end{equation}
(we agree that $x_H^\infty = 0$).
\end{definition}

Finally, the sections of our quasi-homogeneous Lie algebroids will
consist of vector fields that are ``adapted'' to each hyperface, but
the data at each hyperface must be compatible.

\begin{definition}\
A {\em boundary Lie datum} ${\mathcal D} = (E_i^H, \phi_H,d_H)$ on a
manifold with corners $M$, where $H$ ranges through the set of
hyperfaces of $M$, consists of:

(i) Lie flags $E_0^H \subset \ldots \subset E_{l_H}^H \subset TH$ such that 
all intersections $E_{i_1}^{H_1}\cap \ldots \cap E_{i_t}^{H_t}$, 
$i_j \in \{ 0, 1, \ldots, l_{H_j}, \infty \}$, as well as 
all finite sums of such intersections, have constant rank on the
set where they are defined (and hence they are vector bundles), and form a distributive
lattice.

(ii) Diffeomorphisms 
$$
	\phi_H =(\pi_H,x_H): V_H \to H \times [0,\epsilon)
$$
such that $\pi_H \pi_{H'} = \pi_{H'} \pi_{H}$ and $x_H \pi_{H'} = x_H$ on
$V_H \cap V_{H'}$, for all hyperfaces $H$ and $H'$.

(iii) Degrees $d_H \in \NN \cup \{\infty\}$.
\end{definition}

If ${\mathcal D} = (E_i^H, \phi_H,d_H)$ is a boundary Lie datum on
$M$, then the diffeomorphism $\phi_H$ defines a complement $NH$
to $TH$ in $TM\vert_H$.

\begin{proposition}\label{prop.algebroid}\ 
Let ${\mathcal D} = (E_i^H, \phi_H,d_H)$ be a boundary Lie datum on
the manifold with corners $M$. Suppose that  $\mathcal F$ is a foliation 
of $M$ such that 
\begin{equation*}
T{\mathcal F}\vert_H = \begin{cases}
E_{l_H}^H + NH, & \text{ if } d_H < \infty,\\
E_{l_H}^H, & \text{ if }   d_H = \infty.
\end{cases}
\end{equation*}  
If we define
$$
\mathcal{A_D}:=\{ X \in \Gamma(T{\mathcal F}), \; X \text{ is }
(E_i^H, \phi_H,d_H)\text{--adapted,  for each hyperface } H\}
$$
and $[x_j\pa_{x_j},\mathcal{A_D}] \subset \mathcal{A_D}$, then
$\mathcal{A_D}$ is a Lie algebra and a projective $\CI(M)$--module.
\end{proposition}

\begin{proof}\ 
The set of points satisfying \eqref{eq.adapted} at a face $H$ is
closed under the Lie bracket, by definition of a Lie flag. Since
$\Gamma(T{\mathcal F})$ is also closed under the Lie bracket, it
follows that $\mathcal{A_D}$ is a Lie algebra.

Fix a point $x_0$ in the interior of a face $F$ of
codimension $k$ with defining functions $x_1,\ldots,x_k$. Let $H_1,
\ldots, H_k$ be the corresponding hyperfaces, $H_j=\{x_j=0\}$, and let $d_i
= d_{H_i}$. Also let
$$
	\pi=\pi_{H_1} \ldots \pi_{H_k} : V:=\cap V_{H_j} \rightarrow F,
$$
where the order of composition is not important because $\pi_H \pi_{H'}=
\pi_{H'}\pi_H$, by the definition of boundary Lie datum.

Let $\nu = (\nu_1, \ldots, \nu_k)$, $0 \le \nu_i \le l_{H_i}$, be a
multi-index and 
$$
	E_\nu = \pi^*(E_{\nu_1}^{H_1} \cap \ldots \cap
	E_{\nu_k}^{H_k}),
$$ 
which, we recall, is a vector bundle on $H_1 \cap
\ldots \cap H_k$, again by the definition of a boundary Lie datum.  Also,
let 
$$
	\nu^{(i)} = (\nu_1, \ldots, \nu_{i-1}, \nu_i -1, \nu_{i+1},
	\ldots , \nu_k)
$$
and choose, arbitrarily, a complement $Y_\nu$ to $\sum E_{\nu^{(i)}}$
in $E_\nu$. 

Let $\mathcal Z$ be the set of vector fields $\{x_i^{d_{i}+1}\pa_{x_i}, 
d_{i}< \infty\}$
defined using the diffeomorphism 
$(\pi, x_1, \ldots, x_k) : V = \cap V_{H_i} \to F \times [0,\epsilon)^k,$
for some small $\epsilon$.
Then the restriction of $\mathcal{A_D}$ to $V$ is 
\begin{multline}\label{eq.onV}
	\CI(V)\mathcal{A_D}\vert_V= \left (\sum_{\nu=(\nu_1, \ldots,
	\nu_k)} x_1^{\nu_1} \ldots x_k^{\nu_k} \Gamma(\pi^*(Y_{\nu}))
	\right ) \bigoplus \left (\bigoplus_{Z_i \in
	\mathcal{Z}}\CI(V)Z_i \right ) \\ \simeq \Gamma(\oplus_\nu
	Y_\nu \oplus \RR^{\mathcal{Z}}).
\end{multline}
Since this is a $\CI(V)$-projective module--the module of sections of
a bundle isomorphic to the direct sum of $\oplus_{\nu}Y_{\nu}$ and the
trivial bundle generated by the set $x_i^{d_i+1}\pa_{x_i}$ ($d_i <
\infty$), we obtain that $\mathcal{A_D}$ is a projective
$\CI(M)$--module, as desired.
\end{proof}

Let us examine now a particular case of this construction. Assume that
in the definition above ${\mathcal F}=M$, that is, that there exists a
single leaf, and that in the boundary Lie data $d_H=0$ and
$E^H_1=TH$. The only choice then is that of the integrable bundles
$E^H_0$, because the choice of the diffeomorphisms $\phi_H$ is not
important. The conditions that these bundles have to satisfy are that,
for each face $F$ (contained in $k$ distinct hyperplanes $H$, where
$k$ is the codimension of $F$), there exist $2k$ commuting surjective
projections
\begin{equation*}
	p_{FH} : TM\vert_F \rightarrow E^H_0\vert_F, \,\quad \text{ and }\,
	q_{FH} : TM\vert_F \rightarrow TH\vert_F.
\end{equation*}

We denote by $A_{\mathcal D}$ the vector bundle on $M$ with 
sections $\mathcal{A_D}$, defined by the above proposition. 
Because $\mathcal{A_D}$ is a Lie algebra, $A_{\mathcal D}$ is
a Lie algebroid. A Lie algebroid $A$ is called {\em quasi-homogeneous}
if it is isomorphic to one of the form $A_{\mathcal D}$ obtained
as above.

\begin{proposition}\ \label{Prop.q.hom}
Let $A \subset TM$ be a quasi-homogeneous algebroid. Then the set of
{\em open} faces of $M$ defines an $A$-invariant regular
stratification of $M$.  Moreover, if $S = \operatorname{Int}(F)$ is an
open face of $M$, then the Lie algebroid $q_S: A_S := A\vert_S \to TS$
is integrable.
\end{proposition}

\begin{proof} Let ${\mathcal F}$ and $(E_i^H, \phi_H,d_H)$
be the foliation and, respectively, the boundary Lie data defining
$A$.  Observe that the sections of $A$ are vector fields that are
tangent to all faces of $M$, and hence each open face of $M$ is
$A$-invariant. This means that the stratification of $M$ by open faces
is an $A$-invariant stratification of $M$.

Fix a face $F$ of codimension $k$ with defining functions $x_1,\ldots,x_k$,
such that $H_j=\{x_j=0\}$. Let $S :=  \operatorname{Int}(F)$ be the
interior of $F$, and denote $A_S:=A\vert_S$. 

For $\nu_0=(0,\ldots,0)$, the intersection 
$$
	E_{\nu_0}:=E_0^{H_1} \cap \ldots \cap E_0^{H_k} \subset TS
$$ 
is an integrable sub-bundle of $TS$ such that $q_S :A_S \to TS$ maps
$\Gamma(A_S)$ surjectively onto $\Gamma(E_{\nu_0})$.  Moreover, the
map $\pi=\pi_{H_1} \circ \ldots \circ \pi_{H_k}$ defines a splitting
of $A_S$, that is, a Lie algebra morphism $\Gamma(E_{\nu_0}) \to
\Gamma(A_S)$. (This splitting is implicit in the description of
$\mathcal{A_D}$ given in \eqref{eq.onV}.) Consequently, $A_S$ is
integrable, by Proposition \ref{Prop.integrable}.
\end{proof}

\begin{theorem} 
Every quasi-homogeneous algebroid is integrable.
\end{theorem}

\begin{proof} 
This follows from Proposition \ref{Prop.q.hom}
and Theorem \ref{Theorem.application}.  
\end{proof}

Because the construction of the differentiable groupoid integrating a
regular, semi-direct product Lie algebroid depends on the construction
of a differentiable groupoid $\mathcal K$ integrating $\ker (q_S)$, it
is useful to have an explicit construction of $\mathcal K$ in the
above proposition.  To this end, we use the notation in the proof of
Proposition \ref{prop.algebroid}. Let $Z_0 \subset Z :=\{x_i^{d_i +1}
\pa_{x_i}\}$ be the set of normal vectors (with respect to the
decomposition induced by $(\pi,x_1,\ldots,x_k)$) such that
$d_i=d_{H_i}=0$. Also, let $A_0$ be the sub-bundle of $A\vert_S$
generated by $Z_0$ and $A_1$ be the vector bundle generated by the
complement of $Z_0$ in $Z$ and by $E_\nu$, for $\nu\not =
(0,\ldots,0)$.  We then obtain the split exact sequence
$$
0 \lra A_1 \lra \ker(q_S) \lra A_0 \lra 0.
$$
Now, each of the bundles $A_0$ and $A_1$ is a Lie algebroid with
vanishing anchor map. In fact, each of $A_0$ and $A_1$ are bundles of
commutative Lie algebras, and hence are integrable: to integrate them,
we just consider each $A_0$ and $A_1$ as a bundle of commutative Lie
groups. Moreover, $A_1$ is an ideal (in Lie algebra sense) of $\ker
(q_S)$, and the sub-bundle $A_1$ acts by derivations on $A_0$ (the
weights are exactly given by the exponents $\nu=(\nu_{i})$). This
action by derivations of each of the fibers $(A_1)_x$ on $(A_0)_x$
exponentiates to an action of the Lie group with Lie algebra $(A_1)_x$
on the Lie group with Lie algebra $(A_0)_x$. Denote the resulting
semidirect product by ${\mathcal K}_x$.  In this way we obtain on
${\mathcal K}={\mathcal K}_x$ a Lie algebra bundle structure, which is
isomorphic to $\ker (q_S)$, as a fiber bundle, via the exponential
map.

If $A$ is an integrable Lie algebroid, then we denote by
$\GR_A$ a $d$-simply-connected differential groupoid
that integrates $A$, which is unique up to isomorphism. If $A$ is a 
quasi-homogeneous algebroid and $\GR$ is any differential groupoid on 
$M$ with Lie algebroid $A$, then we 
call $\GR$ a {\em quasi-homogeneous} groupoid.

We now make some elementary remarks on quasi-homogeneous algebroids
and groupoids.  If $A \subset M$ is a quasi-homogeneous algebroid and
$F\subset M$ is a face of $M$, then $A_F$ is not a quasi-homogeneous
algebroid unless $F=M$. Let $N_F \subset TM$ be the normal bundle to
$F$, with the trivialization given by the boundary Lie data used in
the definition of $A$. Then $A_F$ is the semidirect product of $A_F \cap TF$ and
$N$. The same will be true of the integrating groupoids. Thus, assume
that $F$ has codimension $k$. Then there exists an action of $\RR^k$
on $\GR_{A_F \cap TF}$, which fixes $F$, the set of units of $\GR_{A_F
\cap TF}$, such that
\begin{equation}
\GR_{A_F} \simeq \GR_{A_F \cap TF} \rtimes \RR^k.
\end{equation}
This relates the differential groupoids $\GR_{A_F}$, associated to the
faces of $M$, to the quasi-homogeneous differential groupoids
$\GR_{A_F \cap TF}$.

Let us now take a closer look at the simplest example of a
quasi-homogeneous, non-regular algebroid, the algebroid ${\mathcal
V}_b(M)$ of all vector fields tangent to the boundary $\pa M \not =
\emptyset$ of a manifold with boundary $M$. The two strata of $M$ are
$\pa M$ and $\operatorname{Int}({M})=M \setminus \pa M$.  If $X$ is a
topological space, we denote by ${\mathcal P}_X$ its {\em path
groupoid}. Then the $d$-simply-connected differential groupoids
$\GR_1$ and $\GR_2$ that integrate ${\mathcal V}_b(M)\vert_{\pa M} =
T(\pa M) \oplus N_{\pa M}$ and ${\mathcal
V}_b(M)\vert_{\operatorname{Int}({M})}$ are, up to isomorphism,
$$
	\GR_1 \simeq {\mathcal P}_{\pa M}\times \RR \, , 
	\ \ \ \text{and} \ \ \
	\GR_2 \simeq {\mathcal P}_{\operatorname{Int}({M})}.
$$

By Theorem \ref{Theorem.Glueying}, or directly, the two groupoids
$\GR_1$ and $\GR_2$ above can be smoothly glued to form a
differentiable groupoid
$$
	\GR_{{\mathcal V}_b(M)} = \GR_1 \cup \GR_2
$$ 
that integrates ${\mathcal V}_b(M)$. This groupoid will be Hausdorff
if, and only if, the morphism $\pi_1(\pa M) \to \pi_1(M)$ is injective. 
The domain map
$$
	d:\GR_{{\mathcal V}_b(M)}  \to M
$$
is a fibration (that is, it has the homotopy lifting property) if, and
only if, $\pi_1(\pa M) \to \pi_1(M)$ is surjective. 

In general, the groupoid $\GR_{{\mathcal V}_b(M)}$ is much larger than
the ``stretched product'' $M^2_b$ considered by Melrose
\cite{Melroseb}, which also gives a groupoid that integrates
${\mathcal V}_b(M)$ after we remove its off-diagonal faces. We get the
same groupoid only if both $\pi_1(\pa M)$ and $\pi_1(M)$ are
trivial. Unlike our $d$-simply connected groupoid $\GR_{{\mathcal
V}_b(M)}$, the groupoid obtained from $M^2_b$ is always Hausdorff and
the domain map $d$ is a fibration.

Let $r \geq 2$ be an integer. Then the same discussion applies to
${\mathcal V}_{b,k}(M)$, the Lie algebra of vector fields that at the
boundary are of the form $x^{r}\pa_x + \sum_j \pa_{y_j}$, for a suitable
coordinate systems $(x,y_1,\ldots,y_{n-1})$ in a neighborhood of a
point of the boundary $\pa M =\{ x= 0\}$ and a suitable choice of a complement
to $T(\pa M)$ in $TM\vert_{\pa M}$.

As a final remark, we now use the algebroid ${\mathcal V}_b(M)$ to
show that some conditions on the groupoids $\GR_S$ are necessary in
Theorem \ref{Theorem.Glueying}. Otherwise the glued groupoid might not
be a smooth manifold. Indeed, consider the same groupoid $\GR_2 \simeq
{\mathcal P}_{\operatorname{Int}({M})}$ for the big-open strata
$\operatorname{Int}({M})$, but a smaller one for the boundary:
$$
\GR_1 \simeq (\pa M \times \pa M) \times \RR.
$$
Then $\GR =\GR_1 \cup \GR_2$ is not a smooth
manifold if $\pi_1(M)$ is non-trivial.

\end{document}